\documentclass[11pt,letterpaper,reqno]{amsart}
\usepackage{amssymb,amscd,amsbsy,mathrsfs}
\newcommand{\lb}{\linebreak}

\renewcommand{\a}{\alpha}
\renewcommand{\b}{\beta}
\newcommand{\g}{\gamma}

\newcommand{\e}{\varepsilon}

\newcommand{\s}{\sigma}

\newcommand{\f}{\varphi}

\newcommand{\B}{{\mathscr B}}

\newcommand{\E}{{\mathscr E}}

\newcommand{\F}{{\mathscr F}}

\newcommand{\X}{{\mathscr X}}
\newcommand{\Y}{{\mathscr Y}}

\newcommand{\R}{{\Bbb R}}
\newcommand{\Z}{{\Bbb Z}}

\newcommand{\bS}{{\boldsymbol S}}

\newcommand{\rf}[1]{(\ref{#1})}

\newcommand{\df}{\stackrel{\mathrm{def}}{=}}

\newcommand{\supp}{\operatorname{supp}}

\newcommand{\trace}{\operatorname{trace}}
\newcommand{\rank}{\operatorname{rank}}
\newcommand{\const}{\operatorname{const}}

\newcommand{\eeq}{\end{equation}}
\newcommand{\beq}{\begin{equation}}
\newcommand{\bay}{\begin{eqnarray}}
\newcommand{\ba}{\begin{align*}}
\newcommand{\ea}{\end{align*}}
\newcommand{\ey}{\end{eqnarray}}
\newcommand{\bey}{\begin{eqnarray*}}
\newcommand{\eey}{\end{eqnarray*}}

\newcommand{\be}{\infty}

\newcommand{\bl}{\blacksquare}

\newcommand{\Pf}{{\bf Proof. }}

\newtheorem{thm}{\hspace{\parindent}Theorem}[section]

\newtheorem{cor}[thm]{\hspace{\parindent}Corollary}
\newtheorem{lem}[thm]{\hspace{\parindent}Lemma}

\usepackage{amsmath}
\usepackage{amsfonts}
\usepackage{amssymb}
\usepackage{graphicx}

\usepackage[T2A,T1]{fontenc}
\DeclareSymbolFont{cyrillic}{T2A}{cmr}{m}{it}
\def\makecyrsymbol#1#2{%
    \begingroup\edef\temp{\endgroup
        \noexpand\DeclareMathSymbol{\noexpand#1}
        {\noexpand\mathalpha}{cyrillic}%
        {\expandafter\expandafter\expandafter
            \calccyr\expandafter\meaning\csname T2A\string#2\endcsname\end}}%
    \temp}
\expandafter\def\expandafter\calccyr\string\char#1\end{#1}
\makecyrsymbol\Zhe\CYRZH
\makecyrsymbol\Be\CYRB
\makecyrsymbol\Shcha\CYRSHCH
\makecyrsymbol\Sha\CYRSH
\makeatletter
\def\upintkern@{\mkern-7mu\mathchoice{\mkern-3.5mu}{}{}{}}
\def\upintdots@{\mathchoice{\mkern-4mu\@cdots\mkern-4mu}%
 {{\cdotp}\mkern1.5mu{\cdotp}\mkern1.5mu{\cdotp}}%
 {{\cdotp}\mkern1mu{\cdotp}\mkern1mu{\cdotp}}%
 {{\cdotp}\mkern1mu{\cdotp}\mkern1mu{\cdotp}}}

\newcommand{\UpMultiIntegral}[1]{%
  \edef\ints@c{\noexpand\upintop
    \ifnum#1=\z@\noexpand\upintdots@\else\noexpand\upintkern@\fi
    \ifnum#1>\tw@\noexpand\upintop\noexpand\upintkern@\fi
    \ifnum#1>\thr@@\noexpand\upintop\noexpand\upintkern@\fi
    \noexpand\upintop
    \noexpand\ilimits@
  }%
  \futurelet\@let@token\ints@a
}
\makeatother

\DeclareFontFamily{OMX}{mdbch}{}
\DeclareFontShape{OMX}{mdbch}{m}{n}{ <->s * [0.8]  mdbchr7v }{}
\DeclareFontShape{OMX}{mdbch}{b}{n}{ <->s * [0.8]  mdbchb7v }{}
\DeclareFontShape{OMX}{mdbch}{bx}{n}{<->ssub * mdbch/b/n}{}

\DeclareSymbolFont{uplargesymbols}{OMX}{mdbch}{m}{n}
\SetSymbolFont{uplargesymbols}{bold}{OMX}{mdbch}{b}{n}
\DeclareMathSymbol{\upintop}{\mathop}{uplargesymbols}{82}
\DeclareMathSymbol{\upointop}{\mathop}{uplargesymbols}{"48}

\DeclareFontEncoding{MDB}{}{}
\DeclareFontFamily{MDB}{mdbch}{}
\DeclareFontShape{MDB}{mdbch}{m}{n}{ <->s * [0.8]  mdbchrmb }{}
\DeclareFontShape{MDB}{mdbch}{b}{n}{ <->s * [0.8]  mdbchbmb }{}
\DeclareFontShape{MDB}{mdbch}{bx}{n}{<->ssub * mdbch/b/n}{}
\DeclareFontSubstitution{MDB}{cmr}{m}{n}
\DeclareSymbolFont{mathdesignB}{MDB}{mdbch}{m}{n}%
\SetSymbolFont{mathdesignB}{bold}{MDB}{mdbch}{b}{n}%
\DeclareMathSymbol{\upintclockwise}{\mathop}{mathdesignB}{128}
\DeclareMathSymbol{\upointclockwise}{\mathop}{mathdesignB}{130}
\DeclareMathSymbol{\upointctrclockwise}{\mathop}{mathdesignB}{132}
\DeclareMathSymbol{\upoiint}{\mathop}{mathdesignB}{134}
\DeclareMathSymbol{\upoiiint}{\mathop}{mathdesignB}{136}

\makeatletter
\newcommand{\upint}{\DOTSI\upintop\ilimits@}
\newcommand{\upoint}{\DOTSI\upointop\ilimits@}
\makeatother

\theoremstyle{remark}

\newtheorem*{rem*}{Remark}

\renewcommand{\f}{{\varphi}}
\renewcommand{\b}{{\beta}}

\newcommand\dg{\frak D}

\newcommand\mB{\mathcal{B}}

\newcommand{\ri}{{\rm i}}

\newcommand{\Bs}{\Be_{\be,1}^1}

\begin{document}

%$\Bs$
%
%%Let $a=\bes =\Zhe=\Sha$, $\Shcha$
%\[
%\Be^p_{1/p}, \quad \| f\|_{\Be^p_{1/p}}
%\]
%
%$\Zhe$, $\Shcha$, $\Sha$

%\end{document}

%\theoremstyle{remark}
%\newtheorem{rem}[thm]{Remark}
%\newtheorem*{rem*}{Remark}

\numberwithin{equation}{section}

\numberwithin{equation}{section}

\title{Functions of perturbed noncommuting\\ unbounded self-adjoint operators}
\author{A.B. Aleksandrov and V.V. Peller}
\thanks{The research is supported by a grant of the Government of the Russian Federation for the state support of scientific research, carried out under the supervision of leading scientists, agreement  075-15-2021-602.
The research is also supported by RFBR [grant number 20-01-00209a] }
\thanks{Corresponding author: V.V. Peller; email: peller@math.msu.edu}

%\date{\today}

\

\begin{abstract}
Let $f$ be a function on $\R^2$ in the inhomogeneous Besov space $\Bs(\R^2)$.
For a pair $(A,B)$ of not necessarily bounded and not necessarily commuting self-adjoint operators, we define the function $f(A,B)$ of $A$ and $B$ as a densely defined linear operator. We show that if $1\le p\le2$, $(A_1,B_1)$ and $(A_2,B_2)$ are pairs of not necessarily bounded and not necessarily commuting self-adjoint operators such that both
$A_1-A_2$ and $B_1-B_2$ belong to the Schatten--von Neumann class $\bS_p$
and $f\in\Bs(\R^2)$, then the following Lipschitz type estimate holds:
$$
\|f(A_1,B_1)-f(A_2,B_2)\|_{\bS_p}
\le\const\|f\|_{\Bs}\max\big\{\|A_1-A_2\|_{\bS_p},\|B_1-B_2\|_{\bS_p}\big\}.
$$
\end{abstract} 

\maketitle

%{\bf
%\footnotesize
%\tableofcontents
%\normalsize
%}

\setcounter{section}{0}
\section{\bf Introduction}
\setcounter{equation}{0}
\label{In}

\

In this paper we obtain Lipschitz type estimates for functions of pairs of unbounded noncommuting self-adjoint operators in the Schatten--von Neumann norms $\bS_p$, $1\le p\le2$. This will
extend the results of \cite{ANP} for functions of pairs of bounded noncommuting self-adjoint operators. 

Recall (see \cite{ANP}) that for a pair $(A,B)$ of not necessarily bounded self-adjoint operators on Hilbert space and for a complex function $f$ on $\R\times\R$ that is a {\it Schur multiplier} with respect to arbitrary spectral measures on $\R$, the function $f(A,B)$ of not necessarily commuting 
self-adjoint operators $A$ and $B$ is defined as the double operator integral 
\bay
\label{fAB}
f(A,B)\df\iint_{\R\times\R}f(x,y)\,dE_A(x)\,dE_B(y)
\df\iint_{\R\times\R}f(x,y)\,dE_A(x)I\,dE_B(y),
\ey
where $I$ is the identity operator. In this case $f(A,B)$ is a bounded linear operator.

Note that double operator integrals
\bay
\label{dvoinoi}
\iint_{\X\times\Y}\Phi(x,y)\,dE_1(x)Q\,dE_2(y)
\ey
appeared first in \cite{DK}. Later Biman and Solomyak created a beautiful theory of double operator integrals, see \cite{BS1}, \cite{BS2} and \cite{BS3}. Here $E_1$ and $E_2$ are spectral measures on Hilbert space, $Q$ is a bounded linear operator on Hilbert space and $\Phi$ is a measurable function of two variables. The double operator integral \rf{dvoinoi} makes sense for an arbitrary bounded linear operator $Q$ if and only if $\Phi$ is a Schur multiplier with respect to the spectral measures $E_1$ and $E_2$. We refer the reader to 
\cite{Pe1} and \cite{AP4} for the definition of Schur multipliers with respect to spectral measures. Note here (see \cite{Pe1} and \cite{AP4}) that $\Phi$ is a Schur multiplier with resect to $E_1$ and $E_2$ if and only if $\Phi$ belongs to the {\it Haagerup tensor product} 
$L^\be(E_1)\!\otimes_{\rm h}\!L^\be(E_2)$, i.e., $\Phi$ admits a representation of the form 
\bay
\label{Phifnpsin}
\Phi(x,y)=\sum_n\f_n(x)\psi_n(y),
\ey
where $\f_n\in L^\be(E_1)$, $\psi_n\in L^\be(E_2)$ and
\bay 
\label{koren'izproizvedeniya}
\left(\left\|\sum_n|\f_n(x)|^2\right\|_{L^\be(E_1)}
\left\|\sum_n|\psi_n(x)|^2\right\|_{L^\be(E_2)}\right)^{1/2}<\be.
\ey
The norm of $\Phi$ in $L^\be(E_1)\!\otimes_{\rm h}\!L^\be(E_2)$ is, by definition, the infimum of \rf{koren'izproizvedeniya} over all representations of $\Phi$ of the form
\rf{Phifnpsin}. For such functions $\Phi$ the double operator integral \rf{dvoinoi} is equal to
$$
\sum_n\left(\int\f_n\,dE_1\right)Q\left(\int\psi_n\,dE_2\right)
$$
and the series converges in the weak operator topology. Moreover, the following estimate holds:
$$
\left\|\iint\Phi\,dE_1Q\,dE_2\right\|\le\|\Phi\|_{L^\be\!\otimes_{\rm h}\!L^\be}\|Q\|
$$
(see \cite{AP4}).

This allows us to define functions $f(A,B)$ of not necessarily commuting bounded self-adjoint operators $A$ and $B$ for functions $f$ on $\R^2$ of class 
$\B^\be\!\otimes_{\rm h}\!\B^\be$. Here $\B^\be$ is the space of bounded Borel functions on $\R$ and the Haagerup tensor product $\B^\be\!\otimes_{\rm h}\!\B^\be$ can be defined by analogy with the definition of $L^\be(E_1)\!\otimes_{\rm h}\!L^\be(E_2)$. 
For such functions $f$, the operator $f(A,B)$ is bounded and
$$
\|f(A,B)\|\le\|f\|_{\B^\be\!\otimes_{\rm h}\!\B^\be}\quad\mbox{for}\quad 
f\in\B^\be\!\otimes_{\rm h}\!\B^\be.
$$

In this paper under a certain assumption on $f$, we define functions $f(A,B)$ for not necessarily bounded and not necessarily commuting self-adjoint operators $A$ and $B$. In this case $f(A,B)$ is a not necessarily bounded densely defined operator.

Recall that it was shown in \cite{ANP} that for pairs 
$(A_1,B_1)$ and $(A_2,B_2)$ of not necessarily commuting bounded self-adjoint operators $A$ and $B$ and for a function $f$ in the {\it homogeneous} Besov space
$B_{\be,1}^1(\R^2)$, the functions of operators $f(A_1,B_1)$ and $f(A_2,B_2)$ are defined and the following Lipschitz type estimate holds in the Schatten-von Neumann spaces $\bS_p$ for $1\le p\le2$:
\bay
\label{vernomezhdu1i2}
\|f(A_1,B_1)-f(A_2,B_2)\|_{\bS_p}\le\const\|f\|_{B_{\be,1}^1}
\max\big\{\|A_1-B_1\|_{\bS_p},\|A_2-B_2\|_{\bS_p}\big\}.
\ey
It was also shown in \cite{ANP} that this inequality is false for $p>2$ and is false in the operator norm.

Recall also that for functions of single self-adjoint operators such such Lipschitz type estimates hold for $1\le p\le\be$, see \cite{Pe1} and \cite{Pe2}. H\"older type estimates for functions of single self-adjoint operators were obtained in \cite{AP1} and \cite{AP2}.

The main purpose of this paper is to establish this inequality for pairs of unbounded self-adjoint operators and for functions $f$ in the {\it inhomogeneous} Besov space $\Bs(\R^2)$.

As in the paper \cite{ANP}, an important role will be played by triple operator integrals. We give a brief introduction to triple operator integrals in \S\,\ref{TrOi}. In \S\,\ref{nekommu} we define functions of pairs of not necessarily bounded and not necessarily commuting operators. We obtain Lipschitz type estmates in the norm of $\bS_p$, $1\le p\le2$, for functions of such pairs  in \S\,\ref{OtsLiptipa}.

Finally, we give in \S\,\ref{Besovy} a brief introduction to Besov spaces.

\

\section{\bf Besov spaces}
\setcounter{equation}{0}
\label{Besovy}

\

In this paper we deal with the homogeneous Besov class $B_{\be,1}^1(\R^2)$ and the inhomogeneous Besov class $\Bs(\R^2)$.  
We refer the reader to the book \cite{Pee} and the papers \cite{ANP} and \cite{AP4} for more information on Besov classes $B_{p,q}^s(\R^d)$. Here we give the definition only in the case
when $p=\be$, $q=s=1$ and $d=2$.
Let $w$ be an infinitely differentiable function on $\R$ such
that
\bay
\label{w}
w\ge0,\quad\supp w\subset\left[\frac12,2\right],\quad\mbox{and} \quad w(t)=1-w\left(\frac t2\right)\quad\mbox{for}\quad t\in[1,2].
\ey

Consider the functions $W_n$, $n\in\Z$, on $\R^2$ such that 
$$
\big(\F W_n\big)(x)=w\left(\frac{\|x\|_2}{2^n}\right),\quad n\in\Z, \quad x=(x_1,x_2),
\quad\|x\|_2\df(x_1^2+x_2^2)^{1/2},
$$
where $\F$ is the {\it Fourier transform} defined on $L^1\big(\R^2\big)$ by
$$
\big(\F f\big)(t)=\!\int\limits_{\R^2} f(x)e^{-{\rm i}(x,t)}\,dx,\!\quad 
x=(x_1,x_2),
\quad t=(t_1,t_2), \!\quad(x,t)\df  x_1t_1+x_2t_2.
$$
Clearly,
$$
\sum_{n\in\Z}(\F W_n)(t)=1,\quad t\in\R^2\setminus\{0\}.
$$

With each tempered distribution $f\in{\mathscr S}^\prime\big(\R^2\big)$, we
associate the sequence $\{f_n\}_{n\in\Z}$,
\bay
\label{fn}
f_n\df f*W_n.
\ey
The formal series
$
\sum_{n\in\Z}f_n
$
is a Littlewood--Paley type expansion of $f$. 
This series does not necessarily converge to $f$. 
%For functions $f$ in $B_{\be,1}^1(\R^d)$,
%\bay
%\label{f(x)-f(y)}
%f(x)-f(y)=\sum_{n\in\Z}\big(f_n(x)-f_n(y)\big),\quad x,~y\in\R^d,
%\ey
%and the series on the right converges uniformly. 

Initially we define the homogeneous Besov class $\dot B^1_{\be,1}\big(\R^2\big)$ as the space of 
$f\in{\mathscr S}^\prime(\R^2)$
such that
\bay
\label{<be}
\|f\|_{B^1_{\be,1}}\df\sum_{n\in\Z}2^n\|f_n\|_{L^\be}<\be.
\ey
According to this definition, the space $\dot B^1_{\be,1}(\R^2)$ contains all polynomials
and all polynomials $f$ satisfy the equality $\|f\|_{B^1_{\be,1}}=0$. Moreover, the distribution $f$ is determined by the sequence $\{f_n\}_{n\in\Z}$
uniquely up to a polynomial. It is easy to see that the series 
$\sum_{n\ge0}f_n$ converges in ${\mathscr S}^\prime(\R^2)$.
However, the series $\sum_{n<0}f_n$ can diverge in general. Obviously, the series
\bay
\label{ryad}
\sum_{n<0}\frac{\partial f_n}{\partial x_1}\quad\mbox{and}\quad
\sum_{n<0}\frac{\partial f_n}{\partial x_2}
\ey
converge uniformly on $\R^2$ for an arbitrary function $f$ in $B^1_{\be,1}(\R^2)$. 
Now we say that  $f$ belongs to the {\it (reduced) homogeneous Besov class} $B^1_{\be,1}(\R^2)$ if \rf{<be} holds and 
$$
\frac{\partial f}{\partial x_j}=\sum_{n\in\Z}\frac{\partial f_n}{\partial x_j},\quad
j=1,\,2.
$$

A function $f$ is determined uniquely by the sequence $\{f_n\}_{n\in\Z}$ up
to a a constant, and a polynomial $g$ belongs to 
$B^1_{\be,1}\big(\R^2\big)$ if and only if it is a constant.

Unlike our previous papers, in this paper we are going to use the inhomogeneous Besov space $\Bs(\R^2)$. We define the function $W^{[0]}$ by
$$
\big(\F W^{[0]}\big)(t)=1-\sum_{n\ge1}(\F W_n)(t),\quad t\in\R^2. 
$$
It is easy to see that
$$
(\F W^{[0]}\big)(t)=1\quad\mbox{if}\quad\|t\|_2\le1\quad\mbox{and}\quad
\supp\F W^{[0]}\subset\big\{t\in\R^2:~\|t\|_2\le2\big\}.
$$
We say that a tempered distribution $f$ on $\R^2$ belongs to the {\it inhomogeneous Besov space} $\Bs(\R^2)$ if
\bay
\label{fnol'}
f^{[0]}\df f*W^{[0]}\in L^\be(\R^2)\quad\mbox{and}\quad
\sum_{n\ge1}2^n\|f_n\|_{L^\be}<\be.
\ey
The norm of a function $f$ in $\Bs(\R^2)$ can be defined by
$$
\|f\|_{\Bs}\df\big\|f^{[0]}\big\|_{L^\be}+\sum_{n\ge1}2^n\|f_n\|_{L^\be}.
$$

Clearly, for $f\in\Bs(\R^2)$, we have
$$
\frac{\partial f}{\partial x_j}=\frac{\partial f^{[0]}}{\partial x_j}
+\sum_{n\ge1}\frac{\partial f_n}{\partial x_j},\quad
j=1,\,2.
$$
Moreover, the series on the right converges uniformly. It is easy to see that $f\in\Bs(\R^2)$ is a dense subset of $B^1_{\be,1}(\R^2)$.

Note that our terminology differs from the standard terminology (see e.g., \cite{Pee}).
However, we believe that our terminology is more convenient.

It is well known and it is easy to see that  $\Bs(\R^2)=B^1_{\be,1}(\R^2)\bigcap L^\be(\R^2)$. 

\

\section{\bf Triple operator integrals and Haagerup-like tensor products}
\setcounter{equation}{0}
\label{TrOi}

\

Triple operator integrals are expressions of the following form:
\bay
\label{troinoi}
\int\limits_{\X_1}\int\limits_{\X_2}\int\limits_{\X_3} 
\Psi(x_1,x_2,x_3)\,dE_1(x_1)T\,dE_2(x_2)R\,dE_3(x_3),
\ey
where $\Psi$ is a bounded measurable function on $\X_1\times\X_2\times\X_3$; 
$E_1$, $E_2$ and $E_3$ are spectral measures on Hilbert space, 
and $T$ and $R$ are bounded linear operators.
Such triple operator integrals can be defined under certain assumptions on $\Psi$, $T$, and $R$.

In \cite{Pe3} triple operator integrals of the form \rf{troinoi} were defined for arbitrary bounded linear operators $T$ and $R$ and for functions $\Psi$ in the integral projective tensor product $L^\be(E_1)\!\otimes_{\rm i}\!L^\be(E_2)\!\otimes_{\rm i}\!L^\be(E_3)$
 of $L^\be(E_1)$, $L^\be(E_2)$ and $L^\be(E_3)$. Moreover, in this case the following inequality holds:
 $$
 \left\|\iiint\Psi\,dE_1T\,dE_2)R\,dE_3\right\|_{\bS_r}\le
 \|\Psi\|_{L^\be\otimes_{\rm i}L^\be\otimes_{\rm i}L^\be}\|T\|_{\bS_p}\|R\|_{\bS_q},\quad
 T\in\bS_p,~R\in\bS_q,
 $$
 $$
 \frac1r=\frac1p+\frac1q\quad\mbox{whenever}\quad\frac1p+\frac1q\le1.
 $$

Later in \cite{JTT} triple operator integrals were defined for functions in the Haagerup tensor product of $L^\be$ spaces.
The {\it Haagerup tensor product}
$L^\be(E_1)\!\otimes_{\rm h}\!L^\be(E_2)\!\otimes_{\rm h}\!L^\be(E_3)$ can be defined as the space of functions $\Psi$ of the form
\bay
\label{htr}
\Psi(x_1,x_2,x_3)=\sum_{j,k\ge0}\a_j(x_1)\b_{jk}(x_2)\g_k(x_3),
\ey
where $\a_j$, $\b_{jk}$ and $\g_k$ are measurable functions such that
\bay
\label{ogr}
\{\a_j\}_{j\ge0}\in L_{E_1}^\be(\ell^2), \quad 
\{\b_{jk}\}_{j,k\ge0}\in L_{E_2}^\be({\mathcal B}),\quad\mbox{and}\quad
\{\g_k\}_{k\ge0}\in L_{E_3}^\be(\ell^2),
\ey
where ${\mathcal B}$ is the space of matrices that induce bounded linear operators on $\ell^2$ and this space is equipped with the operator norm.

The norm of $\Psi$ in 
$L^\be\!\otimes_{\rm h}\!L^\be\!\otimes_{\rm h}\!L^\be$ is, by definition, 
the infimum of
$$
\|\{\a_j\}_{j\ge0}\|_{L^\be(\ell^2)}\|\{\b_{jk}\}_{j,k\ge0}\|_{L^\be({\mB})}
\|\{\g_k\}_{k\ge0}\|_{L^\be(\ell^2)}
$$
over all representations of $\Psi$ of the form \rf{htr} with $\{\a_j\}_{j\ge0}$,
$\{\b_{jk}\}_{j,k\ge0}$ and $\{\g_k\}_{k\ge0}$ satisfying \rf{ogr}.

Let 
$\Psi\in L^\be\!\otimes_{\rm h}\!L^\be\!\otimes_{\rm h}\!L^\be$ and  
suppose that \rf{htr} and \rf{ogr} hold. The triple operator integral \rf{troinoi} is defined by
\begin{align}
\label{htraz}
\iiint\Psi(x_1,x_2,x_3)&\,dE_1(x_1)T\,dE_2(x_2)R\,dE_3(x_3)\nonumber\\[.2cm]
=&
\sum_{j,k\ge0}\left(\int\a_j\,dE_1\right)T\left(\int\b_{jk}\,dE_2\right)
R\left(\int\g_k\,dE_3\right)\nonumber\\[.2cm]
=&\lim_{M,N\to\be}~\sum_{j=0}^N\sum_{k=0}^M
\left(\int\a_j\,dE_1\right)T\left(\int\b_{jk}\,dE_2\right)
R\left(\int\g_k\,dE_3\right).
\end{align}

The series in \rf{htraz} converges in the weak operator topology, the sum of the series does not depend on the choice of a representation
\rf{htr}, it determines a bounded linear operator and
\bay
\label{ogra}
\left\|\iiint\Psi(x_1,x_2,x_3)\,dE_1(x_1)T\,dE_2(x_2)R\,dE_3(x_3)\right\|
\le\|\Psi\|_{L^\be\!\otimes_{\rm h}\!L^\be\!\otimes_{\rm h}\!L^\be}
\|T\|\cdot\|R\|
\ey
(see \cite{ANP} and \cite{AP5}).

It was shown in \cite{ANP} (see also \cite{AP5}) that in the case when 
$$
T\in\bS_p,\quad R\in\bS_q,\quad p,\,q\in[2,\be]\quad\mbox{and}\quad
\Psi\in L^\be\!\otimes_{\rm h}\!L^\be\!\otimes_{\rm h}\!L^\be,
$$
the triple operator integral
\rf{troinoi} belongs to $\bS_r$, $1/r=1/p+1/q$, and
\bay
\label{nerdlyaHaagerupa}
\left\|\iiint\!\Psi(x_1,x_2,x_3)\,dE_1(x_1)T\,dE_2(x_2)R\,dE_3(x_3)\right\|_{\bS_r}
\!\!\!\le\|\Psi\|_{L^\be\!\otimes_{\rm h}\!L^\be\!\otimes_{\rm h}\!L^\be}
\|T\|_{\bS_p}\|R\|_{\bS_q}.
\ey
Note that for $p=\be$, by $\|\cdot\|$ we mean the operator norm.

It turns out that Lipschitz type estimates for functions of pairs of noncommuting operators depend on estimates of triple operator integrals with integrands in {\it Haagerup-like tensor products} of $L^\be$ spaces. Such tensor products were introduced in \cite{ANP}.

\medskip

{\bf Definition 1.} 
{\it A function $\Psi$ of three variables is said to belong to the Haagerup-like tensor product 
$L^\be(E_1)\!\otimes_{\rm h}\!L^\be(E_2)\!\otimes^{\rm h}\!L^\be(E_3)$ of the first kind if it admits a representation
\bay
\label{yaH}
\Psi(x_1,x_2,x_3)=\sum_{j,k\ge0}\a_j(x_1)\b_{k}(x_2)\g_{jk}(x_3),\quad x_j\in\X_j,
\ey
with $\{\a_j\}_{j\ge0},~\{\b_k\}_{k\ge0}\in L^\be(\ell^2)$ and 
$\{\g_{jk}\}_{j,k\ge0}\in L^\be(\mB)$. As usual, 
$$
\|\Psi\|_{L^\be\otimes_{\rm h}\!L^\be\otimes^{\rm h}\!L^\be}
\df\inf\big\|\{\a_j\}_{j\ge0}\big\|_{L^\be(\ell^2)}
\big\|\{\b_k\}_{k\ge0}\big\|_{L^\be(\ell^2)}
\big\|\{\g_{jk}\}_{j,k\ge0}\big\|_{L^\be(\mB)},
$$
the infimum being taken over all representations of the form {\em\rf{yaH}}}.

\medskip

Let us now define triple operator integrals whose integrand belong to the tensor product
$L^\be(E_1)\!\otimes_{\rm h}\!L^\be(E_2)\!\otimes^{\rm h}\!L^\be(E_3)$.

Let $1\le p\le2$. For 
$\Psi\in L^\be(E_1)\!\otimes_{\rm h}\!L^\be(E_2)\!\otimes^{\rm h}\!L^\be(E_3)$, for a bounded linear operator $R$, and for an operator $T$ of class $\bS_p$, we define the triple operator integral
\bay
\label{WHft}
W=\iint\!\!\upint\Psi(x_1,x_2,x_3)\,dE_1(x_1)T\,dE_2(x_2)R\,dE_3(x_3)
\ey
as the following continuous linear functional on $\bS_{p'}$,
$1/p+1/p'=1$ (on the class of compact operators in the case $p=1$):
\bay
\label{fko}
Q\mapsto
\trace\left(\left(
\iiint
\Psi(x_1,x_2,x_3)\,dE_2(x_2)R\,dE_3(x_3)Q\,dE_1(x_1)
\right)T\right).
\ey

\medskip

Clearly, the triple operator integral in \rf{fko} is well defined because the function
$$
(x_2,x_3,x_1)\mapsto\Psi(x_1,x_2,x_3)
$$ 
belongs to the Haagerup tensor product 
$L^\be(E_2)\!\otimes_{\rm h}\!L^\be(E_3)\!\otimes_{\rm h}\!L^\be(E_1)$. It follows easily from \rf{nerdlyaHaagerupa} that
\bay
\label{WsSp}
\|W\|_{\bS_p}\le\|\Psi\|_{L^\be\otimes_{\rm h}\!L^\be\otimes^{\rm h}\!L^\be}
\|T\|_{\bS_p}\|R\|,\quad1\le p\le2.
\ey

%It is easy to see that in the case when $\Psi$ belongs to the projective tensor product $L^\be(E_1)\hat\otimes L^\be(E_2)\hat\otimes L^\be(E_3)$, the definition of the triple operator integral given above is consistent with the definition of the triple operator integral given in \rf{otoi}. Indeed, it suffices to verify this for functions $\Psi$ of the form
%$$
%\Psi(x_1,x_2,x_3)=\f(x_1)\psi(x_2)\chi(x_3),\quad\f\in L^\be(E_1),\quad
%\psi\in L^\be(E_2),\quad\chi\in L^\be(E_3),
%$$
%in which case the verification is obvious.

We also need triple operator integrals in the case when $T$ is a bounded linear operator and $R\in\bS_p$, $1\le p\le2$.

\medskip

{\bf Definition 2.} {\it A function $\Psi$ is said to belong to the Haagerup-like tensor product $L^\be(E_1)\!\otimes^{\rm h}\!L^\be(E_2)\!\otimes_{\rm h}\!L^\be(E_3)$
of the second kind if
$\Psi$ admits a representation
\bay
\label{preds}
\Psi(x_1,x_2,x_3)=\sum_{j,k\ge0}\a_{jk}(x_1)\b_{j}(x_2)\g_k(x_3)
\ey
where $\{\b_j\}_{j\ge0},~\{\g_k\}_{k\ge0}\in L^\be(\ell^2)$, 
$\{\a_{jk}\}_{j,k\ge0}\in L^\be(\mB)$. The norm of $\Psi$ in 
the space $L^\be\otimes^{\rm h}\!L^\be\otimes_{\rm h}\!L^\be$ is defined by
$$
\|\Psi\|_{L^\be\otimes^{\rm h}\!L^\be\otimes_{\rm h}\!L^\be}
\df\inf\big\|\{\a_j\}_{j\ge0}\big\|_{L^\be(\ell^2)}
\big\|\{\b_k\}_{k\ge0}\big\|_{L^\be(\ell^2)}
\big\|\{\g_{jk}\}_{j,k\ge0}\big\|_{L^\be(\mB)},
$$
the infimum being taken over all representations of the form {\em\rf{preds}}}.

\medskip

Suppose now that 
$\Psi\in L^\be(E_1)\!\otimes^{\rm h}\!L^\be(E_2)\!\otimes_{\rm h}\!L^\be(E_3)$,
$T$ is a bounded linear operator, and $R\in\bS_p$, $1\le p\le2$. The continuous linear functional 
$$
Q\mapsto
\trace\left(\left(
\iiint\Psi(x_1,x_2,x_3)\,dE_3(x_3)Q\,dE_1(x_1)T\,dE_2(x_2)
\right)R\right)
$$
on the class $\bS_{p'}$ (on the class of compact operators in the case $p=1$) 
determines an operator $W$ of class $\bS_p$, which
we call the triple operator integral
\bay
\label{WHst}
W=\upint\!\!\!\iint\Psi(x_1,x_2,x_3)\,dE_1(x_1)T\,dE_2(x_2)R\,dE_3(x_3).
\ey

Moreover,
\bay
\label{eshchorazW}
\|W\|_{\bS_p}\le
\|\Psi\|_{L^\be\otimes^{\rm h}\!L^\be\otimes_{\rm h}\!L^\be}
\|T\|\cdot\|R\|_{\bS_p}.
\ey

In \cite{ANP} we obtained more general Schatten--von Neumann estimates of triple operator integrals with integrands in the Haagerup-like tensor product of $L^\be$ spaces. 
Later those estimates were generalized in \cite{AP5} to the most general situation.

Note that in the same way we can also define Haagerup-like tensor products 
\lb$\B^\be\!\otimes_{\rm h}\!\B^\be\!\otimes^{\rm h}\!\B^\be$ and 
$\B^\be\!\otimes^{\rm h}\!\B^\be\!\otimes_{\rm h}\!\B^\be$. 

For a continuously differentiable function $f$ on $\R^2$, we define
the divided differences $\dg^{[1]}f$ and $\dg^{[2]}f$ by
$$
\big(\dg^{[1]}f\big)(x_1,x_2,y)\df\frac{f(x_1,y)-f(x_2,y)}{x_1-x_2},
\quad\mbox{and}\quad x_1\ne x_2
$$
and
$$
\big(\dg^{[2]}f\big)(x,y_1,y_2)\df\frac{f(x,y_1)-f(x,y_2)}{y_1-y_2},\quad y_1\ne y_2.
$$
In the case when $x_1=x_2$ or $y_1=y_2$ in the definitions of $\dg^{[1]}f$ and 
$\dg^{[2]}f$ we should replace the divided differences with the corresponding partial derivatives. 

It was shown in \cite{ANP} that if $\s>0$ and $f\in\E_\s^\be(\R^2)$, then 
$\dg^{[1]}f\in\B^\be\!\otimes_{\rm h}\!\B^\be\!\otimes^{\rm h}\!\B^\be$ and
$\dg^{[2]}f\in\B^\be\!\otimes^{\rm h}\!\B^\be\!\otimes_{\rm h}\!\B^\be$; moreover,
\bay
\label{otsenkaD1}
\big\|\dg^{[1]}f\big\|_{\B^\be\otimes_{\rm h}\B^\be\otimes^{\rm h}\B^\be}
\le\const\s\|f\|_{L^\be(\R^2)}
\ey
and
\bay
\label{otsenkaD2}
\big\|\dg^{[2]}f\big\|_{\B^\be\otimes^{\rm h}\B^\be\otimes_{\rm h}\B^\be}
\le\const\s\|f\|_{L^\be(\R^2)}.
\ey
Recall that for $\s>0$
$$
\E_\s^\be(\R^2)\df\Big\{f\in L^\be(\R^2):~\supp\F f\subset\{t\in\R^2:~\|t\|_2\le\s\}\Big\}.
$$

These inequalities imply that if $f\in B_{\be,1}^1(\R^2)$, $\dg^{[1]}f\in\B^\be\!\otimes_{\rm h}\!\B^\be\!\otimes^{\rm h}\!\B^\be$ and
$\dg^{[2]}f\in\B^\be\!\otimes^{\rm h}\!\B^\be\!\otimes_{\rm h}\!\B^\be$; moreover,
\bay
\label{otsenkaD3}
\big\|\dg^{[1]}f\big\|_{\B^\be\otimes_{\rm h}\B^\be\otimes^{\rm h}\B^\be}
\le\const\|f\|_{B_{\be,1}^1}
\ey
and
\bay
\label{otsenkaD4}
\big\|\dg^{[2]}f\big\|_{\B^\be\otimes^{\rm h}\B^\be\otimes_{\rm h}\B^\be}
\le\const\|f\|_{B_{\be,1}^1}.
\ey

Next, it was shown in \cite{ANP} that if $1\le p\le2$, $f\in B_{\be,1}^1(\R^2)$ and $(A_1,B_1)$ and $(A_2,B_2)$ are pairs of bounded self-adjoint operators such that $A_2-A_1\in\bS_p$ and 
$B_2-B_1\in\bS_p$, then
\begin{align*}
f(A_1,B_1)&-f(A_2,B_2)\nonumber\\[.2cm]
&=
\iint\!\!\upint\frac{f(x_1,y)-f(x_2,y)}{x_1-x_2}
\,dE_{A_1}(x_1)(A_1-A_2)\,dE_{A_2}(x_2)\,dE_{B_1}(y),\nonumber\\[.2cm]
&+\upint\!\!\!\iint\frac{f(x,y_1)-f(x,y_2)}{y_1-y_2}
\,dE_{A_2}(x)\,dE_{B_1}(y_1)(B_1-B_2)\,dE_{B_2}(y_2),
\end{align*}
and so the following Lipschitz type estimate holds:
$$
\|f(A_1,B_1)-f(A_2,B_2)\|_{\bS_p}\le\const\|f\|_{B_{\be,1}^1}
\max\big\{\|A_1-A_2\|_{\bS_p},\|B_1-B_2\|_{\bS_p}\big\}.
$$

In this paper we are going to prove the same inequality without the assumption that the operators $A_1$, $B_1$, $A_2$ and $B_2$ are bounded. However, to define the functions $f(A,B)$ for unbounded self-adjoint operators $A$ and $B$, we impose the assumption that $f\in\Bs(\R^2)$.

\

\section{\bf Functions of pairs of unbounded noncommuting self-adjoint operators}
\label{nekommu}

\

Recall that in \S\,\ref{In} we have defined functions of not necessarily commuting self-adjoint operators by formula \rf{dvoinoi} for functions $f$ in the Haagerup tensor product 
$\B^\be\!\otimes_{\rm h}\!\B^\be$ and for self-adjoint operators $A$ and $B$. Herewith
the following estimate holds:
$$
\|f(A,B)\|\le\|f\|_{\B^\be\otimes_{\rm h}\B^\be},\quad f\in\B^\be\otimes_{\rm h}\B^\be.
$$
It is easy to see that if $f\in\B^\be\!\otimes_{\rm h}\!\B^\be$, $g\in \B^\be$ 
and $f_\flat$ is the function defined by
$f_\flat(s,t)\df g(t)f(s,t)$, 
then $f_\flat\in\B^\be\!\otimes_{\rm h}\!\B^\be$
and $f_\flat(A,B)=f(A,B)g(B)$.
This gives us an idea how to define $f(A,B)$ as a densely defined not necessarily bounded operator in the case when $f\not\in\B^\be\!\otimes_{\rm h}\!\B^\be$.

Let $f_\sharp$ be the function defined by
$f_\sharp(s,t)\df(1-\ri t)^{-1}f(s,t)$
and assume that $f_\sharp\in\B^\be\!\otimes_{\rm h}\!\B^\be$. 
We define $f(A,B)$ by
$$
f(A,B)\df f_\sharp(A,B)(I-\ri B)=\left(\,\,\,\iint\limits_{\R\times\R}f_\sharp(s,t)\,dE_A(s)\,dE_B(t)\right)(I-\ri B).
$$
Note that $f(A,B)$ is defined  {\it on the domain $D(B)$ of $B$}. It is a not necessarily bounded  but the operator $f(A,B)(I-\ri B)^{-1}$ must be bounded.

We could try to give another definition of $f(A,B)$ for self-adjoint operators $A$ and $B$. 
Put $f_\natural(s,t)\df(1-\ri s)^{-1}f(s,t)$. If $f_\natural\in\B^\be\!\otimes_{\rm h}\!\B^\be$,
then we could define $f(A,B)$ as follows
$$
f(A,B)\df (I-\ri A)f_\natural(A,B)= (I-\ri A)\left(\,\,\,\iint\limits_{\R\times\R}f_\natural(s,t)\,dE_A(s)\,dE_B(t)\right).
$$
However, with this definition of $f(A,B)$ the domain $D(f(A,B))$ of $f(A,B)$ does not have to be dense.

Corollary 7.3 of \cite{AP6} implies that for $f\in(\E_\s^\be)(\R^2)$ and $\s>0$, 
the function $f_\sharp$ belongs to $\B^\be\!\otimes_{\rm h}\!\B^\be$, and so 
the operator $f_\sharp(A,B)$ is bounded while
$f(A,B)$ is a possibly unbounded operator defined on the domain of $B$.
Moreover, 
\bay
\label{1+sigma}
\|f_\sharp\|_{\B^\be\otimes_{\rm h}\B^\be}\le\const(1+\s)\|f\|_{L^\be(\R^2)},
\ey
see Corollary 7.3 of \cite{AP6}.

\medskip

{\bf Remark.} Let $B_0$ be a self-adjoint operator such that $B_0-B$ bounded.
Then $f(A,B)(I-\ri B_0)^{-1}$ is a bounded operator  for every $f\in(\E_\s^\be)(\R^2)$, where $\s>0$.

\

\Pf It suffices to observe that
\bey
f(A,B)(I-\ri B_0)^{-1}=f_\sharp(A,B)(I-\ri B)(I-\ri B_0)^{-1}\\
=\ri f_\sharp(A,B)(B_0-B)(I-\ri B_0)^{-1}+f_\sharp(A,B). \quad\bl
\eey

%\begin{lem} 
%\label{cor73}
%Let $f\in\E^\be_\s(\R^2)$, where $\s>0$. Then 
%$$
%\vk(x)f(x,y)\in  \B^\be(\R)\otimes_{\rm h}\B^\be(\R)
%\quad\text{and}\quad \vk(y)f(x,y)\in \B^\be(\R)\otimes_{\rm h}\B^\be(\R).
%$$
%Moreover,
%$$
%\|\vk(x)f(x,y)\|_{\B^\be(\R)\otimes_{\rm h}\B^\be(\R)}\le\const(1+\s)\|f\|_{L^\be(\R^2)}
%$$
%and
%$$
%\|\vk(y)f(x,y)\|_{\B^\be(\R)\otimes_{\rm h}\B^\be(\R)}\le\const(1+\s)\|f\|_{L^\be(\R^2)}.
%$$
%\end{lem}

%\Pf A similar statement for $f\in(\E^\be_\s)_+(\R^2)$ can be founded in \cite{AP6?} (Corollary 7.3).
%The lemma can be proved in the same way.  
%
%Besides, the lemma can be reduced to Corollary 7.3 in \cite{AP6?}.
%It suffices to observe that if $f\in\E^\be_\s(\R^2)$, then 
%$e^{{\rm i}\s x}e^{{\rm i}\s y}f(x,y)\in\Big(\E^\be_{(1+\sqrt2)\s}\Big)_+(\R^2)$. $\bl$
%

%\begin{thm} Let $f(x,y)\in\Bs$. Then $f(x,y)(1-\ri y)^{-1}\in \B^\be\otimes_{\rm h}\B^\be$
%and $\|f(x,y)(1-\ri y)^{-1}\|_{\B^\be\otimes_{\rm h}\B^\be}\le\const\|f(x,y)\|_{\Bs}$.
%\end{thm}

\begin{thm} 
\label{th73}
Let $f\in\Bs(\R^2)$. Then $f_\sharp\in \B^\be\otimes_{\rm h}\B^\be$
and $\|f_\sharp\|_{\B^\be\otimes_{\rm h}\B^\be}\le\const\|f\|_{\Bs}$.
\end{thm}

\Pf Let $f\in\Bs(\R^2)$. Then $f=f^{[0]}+\sum\limits_{n=1}^\infty f_n$ (the series converges in $\Bs(\R^2)$), 
where $f^{[0]}$ and $f_n$ denote the same as in Section \ref{Besovy}. It remains to apply inequality \rf{1+sigma}
to $f^{[0]}$ and $f_n$ for every $n\ge1$. $\bl$

\medskip

Theorem \ref{th73} allows us for $f\in\Bs(\R^2)$, to define the operator $f(A,B)$  on the domain of $B$ for arbitrary self-adjoint operators $A$ and $B$.

%\
%
%{\bf Example 1.} Let $Bf=xf$, where $f\in L^2(\R)$
%Put 
%$$
%(Af)(x)=\left(\int_\R(1-\ri t)^{-1}f(t)\,dt\right)(1+\ri x)^{-1}.
%$$
%Let $u_n=\frac1n\chi_{[0,n]}$.  Then $\lim\limits_{n\to\be}\|u_n\|_{L^2(\R)}=0$
%and $A(I-\ri B)u_n=(1+\ri x)^{-1}$ for all $n$.
%
%\
%
%{\bf Example 2.} Let $A$ be an unbounded self-adjoint operator. Then
%there exists a vector $u\in\h$ such that $u\not\in D(A)$.
%Let $B$ be a self-adjoint operator such that $B(\h)=\C h$.
%Then $D((I-\ri A)B)=\ker B$. 

\

\section{\bf Integral formulae for operator differences and Lipschitz type estimates}
\label{OtsLiptipa}

\

In this section we obtain the main results of the paper. We represent operator differences in terms of triple operator integrals and obtain Lipschitz type estimates for functions of pairs of not necessarily bounded and not necessarily commuting operators in the norms of $\bS_p$, $1\le p \le2$.

Throughout this section by a self-adjoint operator we mean a not necessarily bounded self-adjoint operator. By a pair of self-adjoint operators we mean a pair of not necessarily bounded and not necessarily commuting self-adjoint operators. 

\begin{lem}
\label{ssylku3}
Let $A$ be a self-adjoint operator. Then $\lim\limits_{\e\to0}(I-\ri\e A)^{-1}=I$
in the strong operator topology.
\end{lem}
\Pf It suffices to observe that 
$$
(I-\ri\e A)^{-1}=\int(1-\ri\e t)^{-1}\,dE_A(t),
$$
$|(1-\ri\e t)^{-1}|\le1$ and 
$\lim\limits_{\e\to0}(1-\ri\e t)^{-1}=1$ for all $t\in\R$, see, e.g., \cite{BS0}, Ch. 5, Sect. 3,
Th. 2. $\bl$

\begin{lem}
\label{ssylku}
Let $\{X_n\}_{n\ge1}$ be a sequence of bounded linear operators such that 
\lb$\lim\limits_{n\to\be}X_n=0$
in the strong operator topology. Then
$$
\lim_{n\to\be}\|X_nC\|_{\bS_2}=\lim_{n\to\be}\|CX_n\|_{\bS_2}=0
$$
for every $C\in\bS_2$.
\end{lem}

\Pf The lemma is trivial in the case when $\rank C<+\be$.
It remains to observe that the finite rank operators 
are dense in $\bS_2$. $\bl$

\begin{cor}
\label{ssylku1}
Let $\{X_n\}_{n\ge1}$ be a sequence of bounded linear operators such that $\lim\limits_{n\to\be}X_n=0$
in the strong operator topology and let  $\{C_n\}_{n\ge1}$ be a convergent sequence in $\bS_2$. Then
$$
\lim_{n\to\be}\|X_nC_n\|_{\bS_2}=\lim_{n\to\be}\|C_nX_n\|_{\bS_2}=0.
$$
\end{cor}

\begin{cor}
\label{ssylku2}
Let $\{X_n\}_{n\ge1}$ and  $\{Y_n\}_{n\ge1}$ be sequences of bounded linear operators
such that $\lim\limits_{n\to\infty}X_n=X$ and $\lim\limits_{n\to\infty}Y_n=Y$
in the strong operator topology. Then for every $C\in\bS_2$,
$$
\lim_{n\to\be}X_nCY_n=XCY
$$
in the norm of $\bS_2$.
\end{cor}

Let $A$ be a self-adjoint operator. Put $A(\e)\df A(I-\ri\e A)^{-1}$ for $\e>0$.
It is easy to see that $A(\e)$ is a bounded operator for every $\e>0$.

In what follows for possibly unbounded self-adjoint operators $A_1$ and $A_2$, we say that $A_1-A_2\in\bS_p$ if the operator $A_1-A_2$ extends by continuity to an
operator of class $\bS_p$.

\begin{lem}
\label{spnorm}
Let $A_1$ and $A_2$ be self-adjoint operators such that $A_1-A_2\in\bS_2$.
Then 
$$
\lim_{\e\to0}(A_1(\e)-A_2(\e))=A_1-A_2
$$
in the norm of  $\bS_2$.
\end{lem}

Strictly speaking, the left-hand side of \rf{volnost'} is defined on the whole space but formally, the right-hand side in general is defined on a proper subset of the Hilbert space.
When we write this equality, we mean that the right-hand side extends to the whole space and the extension coincides with the left-hand side.

\medskip

\Pf  We have
\bay
\label{volnost'}
A_1(I-\ri\e A_1)^{-1}-A_2(I-\ri\e A_2)^{-1}=(I-\ri\e A_1)^{-1}(A_1-A_2)(I-\ri\e A_1)^{-1}.
\ey
Since $\lim\limits_{\e\to0}(I-\ri\e A)^{-1}=I$ in the strong operator topology by Lemma \ref{ssylku3},
it remains to apply Corollary \ref{ssylku2}. $\bl$

%\begin{lem}
%\label{spnorm}
%Let $A_1$ and $A_2$ be self-adjoint operators such that $A_2-A_1\in\bS_2$.
%Then 
%$$
%\lim_{\e\to0}(A_1(\e)-A_2(\e))=A_1-A_2
%$$
%in the norm of  $\bS_2$.
%\end{lem}
%
%\Pf  We have
%\bey
%A_1(I-\ri\e A_1)^{-1}-A_2(I-\ri\e A_2)^{-1}=(I-\ri\e A_1)^{-1}(A_1-A_2)(I-\ri\e A_1)^{-1}\\
%\eey
%Hence,
%\bey
%\|A_1(\e)-A_2(\e)\|\le\|\ri\e A_1(I-\ri\e A_1)^{-1}\|\cdot\|A_1-A_2\|_{\bS_2}\cdot\|(I-\ri\e A_1)^{-1}\|\\
%+\|A_1-A_2\|_{\bS_2}\cdot\|\ri\e A_2(I-\ri\e A_2)^{-1}\|\le2
%\eey
%Note that $\lim\limits_{\e\to0}(I-\ri\e A)^{-1}=I$ in the strong operator topology.
%It remains to apply Corollary \ref{ssylku2}. $\bl$

\begin{thm}
\label{teor1}
Let $f\in \E^\be_\s(\R^2)$. Suppose that $A_1$, $A_2$ and $B$ are self-adjoint operators such that 
$A_1-A_2\in\bS_2$. Then the following identity holds:
\begin{align*}
f(A_1,B)&-f(A_2,B)\\[.2cm]
&=
\iint\!\!\upint\frac{f(x_1,y)-f(x_2,y)}{x_1-x_2}
\,dE_{A_1}(x_1)(A_1-A_2)\,dE_{A_2}(x_2)\,dE_{B}(y).
\end{align*}
\end{thm}

\Pf Recall that under the hypotheses of the theorem the divided difference $\dg^{[1]}f\in\B^\be\!\otimes_{\rm h}\!\B^\be\!\otimes^{\rm h}\!\B^\be$, see \rf{otsenkaD1}, and so by Definition 1 in \S\,\ref{TrOi}, the triple operator integral on the right is well defined.

To prove the equality, it suffices to verify that
\begin{align*}
(f(A_1,B)&-f(A_2,B))(I-\ri B)^{-1}\\[.2cm]
&\!\!\!\!\!=
\left(\iint\!\!\upint\frac{f(x_1,y)-f(x_2,y)}{x_1-x_2}
\,dE_{A_1}(x_1)(A_1-A_2)\,dE_{A_2}(x_2)\,dE_{B}(y)\right)(I-\ri B)^{-1}.
\end{align*}

Let us observe that Lemma \ref{spnorm} and Theorem 5.1 in \cite{ANP} imply the following equality
\begin{align*}
\lim_{\e\to0}\left(\iint\!\!\upint\frac{f(x_1,y)-f(x_2,y)}{x_1-x_2}
\,dE_{A_1}(x_1)(A_1(\e)-A_2(\e))\,dE_{A_2}(x_2)\,dE_{B}(y)\right)(I-\ri B)^{-1}\\
=\iint\!\!\upint\frac{f(x_1,y)-f(x_2,y)}{x_1-x_2}
\,dE_{A_1}(x_1)(A_1-A_2)\,dE_{A_2}(x_2)(I-\ri B)^{-1}\,dE_{B}(y)
\end{align*}
in the norm of $S_2$.
On the other hand, we have
\begin{multline*}
\lim_{\e\to0}\left(\iint\!\!\upint\!\frac{f(x_1,y)-f(x_2,y)}{x_1-x_2}\!
\,dE_{A_1}(x_1)(A_1(\e)-A_2(\e))\,dE_{A_2}(x_2)\,dE_{B}(y)\right)(I-\ri B)^{-1}\\
=\!\lim_{\e\to0}\!\!\left(\iint\!\!\upint\!\!\frac{f(x_1,y)\!-\!f(x_2,y)}{x_1-x_2}\!\left(\!\frac{x_1}{1\!-\!\ri\e x_1}\!-\!\frac{x_2}{1\!-\!\ri\e x_2}\!\right)\!
dE_{A_1}(x_1)dE_{A_2}(x_2)dE_{B}(y)\!\right)\!\!(I-\ri  B)^{\!-1}\\
=\lim_{\e\to0}\iint\!\!\upint\frac{f(x_1,y)-f(x_2,y)}{(1-\ri \e x_1)(1-\ri \e x_2)(1-\ri y)}
\,dE_{A_1}(x_1)\,dE_{A_2}(x_2)\,dE_{B}(y)\\
=\lim_{\e\to0}\iint\!\!\upint\frac{f(x_1,y)}{(1-\ri \e x_1)(1-\ri \e x_2)(1-\ri y)}
\,dE_{A_1}(x_1)\,dE_{A_2}(x_2)\,dE_{B}(y)\\
-\lim_{\e\to0}\iint\!\!\upint\frac{f(x_2,y)}{(1-\ri \e x_1)(1-\ri \e x_2)(1-\ri y)}
\,dE_{A_1}(x_1)\,dE_{A_2}(x_2)\,dE_{B}\\
=\lim_{\e\to0}\iint \frac{f(x_1,y)}{(1-\ri \e x_1)(1-\ri y)}
\,dE_{A_1}(x_1)(1-\ri \e A_2)^{-1}\,dE_{B}(y)\\
-\lim_{\e\to0}(1-\ri \e A_1)^{-1}\iint\frac{f(x_2,y)}{(1-\ri \e x_2)(1-\ri y)}
\,dE_{A_2}(x_2)\,dE_{B}(y)
=\lim_{\e\to0} X_\e-\lim_{\e\to0} Y_\e\,
\end{multline*}
where
$$
X_\e\df\lim_{\e\to0}(1-\ri \e A_1)^{-1}\iint\frac{f(x_1,y)}{1-\ri y}
\,dE_{A_1}(x_1)(1-\ri \e A_2)^{-1}\,dE_{B}(y)
$$
and
$$
Y_\e\df
\lim_{\e\to0}(1-\ri \e A_1)^{-1}(1-\ri \e A_2)^{-1}\iint \frac{f(x_2,y)}{1-\ri y}
\,dE_{A_2}(x_2)\,dE_{B}(y).
%\df\lim_{\e\to0} X_\e-\lim_{\e\to0} Y_\e.
$$

Clearly, 
$$
\lim_{\e\to0} Y_\e=\iint \frac{f(x_2,y)}{1-\ri y}\,dE_{A_2}(x_2)\,dE_{B}(y)=f(A_2,B)(I-\ri  B)^{-1}
$$
in the strong operator topology.
To compute the limit $\lim\limits_{\e\to0} X_\e$, we observe that by Proposition 3.3 of \cite{BS3},
$$
\lim_{\e\to0}\iint\frac{f(x_1,y)}{1-\ri y}\,dE_{A_1}(x_1)(1-\ri \e A_2)^{-1}\,dE_{B}(y)
=\iint\frac{f(x_1,y)}{1-\ri y}\,dE_{A_1}(x_1)\,dE_{B}(y)
$$
in the strong operator topology. Hence,
\begin{align*}
\lim_{\e\to0} X_\e=\iint\frac{f(x_1,y)}{1-\ri y}\,dE_{A_1}(x_1)\,dE_{B}(y)
=f(A_1,B)(I-\ri B)^{-1}
\end{align*}
in the strong operator topology. Thus, we have proved that
\bey
\lim_{\e\to0}\left(\iint\!\!\upint\frac{f(x_1,y)-f(x_2,y)}{x_1-x_2}
\,dE_{A_1}(x_1)(A_1(\e)-A_2(\e))\,dE_{A_2}(x_2)\,dE_{B}(y)\right)(I-\ri  B)^{-1}\\
=(f(A_1,B)-f(A_2,B))(I-\ri  B)^{-1}
\eey
in the strong operator topology. $\bl$.

\begin{cor} 
\label{cor57}
Let $f\in \E^\be_\s(\R^2)$, $1\le p\le2$. Suppose that $A_1$, $A_2$ and $B$ are self-adjoint operators such that 
$A_2-A_1\in\bS_p$. Then the following inequality holds:
$$
\|f(A_1,B)-f(A_2,B)\|_{\bS_p}\le\const\s\|f\|_{L^\be(\R^2)}\|A_1-A_2\|_{\bS_p}.
$$
\end{cor}

\Pf Applying \rf{WsSp} and \rf{otsenkaD1}, we obtain
\begin{align*}
\|f(A_1,B)-f(A_2,B)\|_{\bS_p}
&\le\big\|\big(\dg^{[1]}f\big)(x_1,x_2,y)\big\|_{\B^\be\otimes_{\rm h}\B^\be\otimes^{\rm h}\B^\be}\|A_1-A_2\|_{\bS_p}\\
&\le\const\s\|f\|_{L^\be(\R^2)}\|A_1-A_2\|_{\bS_p}.\qquad\qquad\qquad\bl
\end{align*}

\begin{cor} 
\label{cor58}
Let $f\in \Bs(\R^2)$, $1\le p\le2$. Suppose that $A_1$, $A_2$ and $B$ are self-adjoint operators such that 
$A_2-A_1\in\bS_p$. Then the following inequality holds:
$$
\|f(A_1,B)-f(A_2,B)\|_{\bS_p}\le\const\|f\|_{\Bs}\|A_1-A_2\|_{\bS_p}.
$$
\end{cor}

\Pf Let $f^{[0]}$ and $f_n$ be the functions defined by \rf{fnol'} and \rf{fn}.
It is easy to see that if $u\in D(B)$, then
$$
f(A_1,B)u-f(A_2,B)u=f^{[0]}(A_1,B)u-f^{[0]}(A_2,B)u+\sum_{n=1}^\infty(f_n(A_1,B)u-f_n(A_2,B)u)
$$
It remains to observe that
$$
\|f^{[0]}(A_1,B)u-f^{[0]}(A_2,B)u\|\le\const\|f^{[0]}\|_{L^\infty}\|A_1-A_2\|\cdot\|u\|
$$
and
$$
\|f_n(A_1,B)u-f_n(A_2,B)u\|\le\const 2^n\|f_n\|_{L^\infty}\|A_1-A_2\|\cdot\|u\|
$$
by Corollary \ref{cor57}. $\bl$

\begin{thm}
\label{teor2}
Let $f\in \E^\be_\s(\R^2)$. Suppose that $A$, $B_1$ and $B_2$ are self-adjoint operators such that 
$B_2-B_1\in\bS_2$. Then the following identity holds:
\begin{align*}
f(A,B_1)&-f(A,B_2)\nonumber\\[.2cm]
&=
\upint\!\!\!\iint\frac{f(x,y_1)-f(x,y_2)}{y_1-y_2}
\,dE_{A}(x)\,dE_{B_1}(y_1)(B_1-B_2)\,dE_{B_2}(y_2).
\end{align*}
\end{thm}

\Pf The proof is similar to the proof of Theorem \ref{teor1}. It follows from inequality \rf{otsenkaD2} that
the right hand side is well defined.

To prove the equality, it suffices to verify that
\begin{align*}
(f(A,B_1)&-f(A,B_2))(I-\ri B_2)^{-1}\nonumber\\[.2cm]
&\!=\!
\left(\upint\!\!\!\iint\frac{f(x,y_1)-f(x,y_2)}{y_1-y_2}
dE_{A}(x)dE_{B_1}(y_1)(B_1-B_2)dE_{B_2}(y_2)\right)(I-{\rm i}B_2)^{-1}.
\end{align*}

Note that Lemma \ref{spnorm} of this paper and Theorem 5.2 in \cite{ANP} imply the following equality
\begin{align*}
\lim_{\e\to0}\left(\upint\!\!\!\iint\frac{f(x,y_1)-f(x,y_2)}{y_1-y_2}
\,dE_{A}(x)\,dE_{B_1}(y_1)(B_1(\e)-B_2(\e))\,dE_{B_2}(y_2)\right)(I-{\rm i}B_2)^{-1}\\
=\left(\upint\!\!\!\iint\frac{f(x,y_1)-f(x,y_2)}{y_1-y_2}
\,dE_{A}(x)\,dE_{B_1}(y_1)(B_1-B_2)\,dE_{B_2}(y_2)\right)(I-{\rm i}B_2)^{-1}
\end{align*}
in the norm of $\bS_2$. On the other hand, we have
\begin{multline*}
\lim_{\e\to0}\left(\upint\!\!\!\iint\frac{f(x,y_1)-f(x,y_2)}{y_1-y_2}\,dE_{A}(x)
\,dE_{B_1}(y_1)(B_1(\e)-B_2(\e))\,dE_{B_2}(y_2)\right)(I-{\rm i}B_2)^{-1}\\
=\!\lim_{\e\to0}\left(\upint\!\!\!\iint\!\frac{f(x,y_1)\!-\!f(x,y_2)}{y_1-y_2}\!
\left(\!\frac{y_1}{1\!-\!{\rm i}\e y_1}\!-\!\frac{y_2}{1\!-\!{\rm i}\e y_2}\!\right)
\!dE_{A}(x)dE_{B_1}(y_1)dE_{B_2}(y_2)\right)\!(I-{\rm i}B_2)^{\!-1}\\
=\lim_{\e\to0}\upint\!\!\!\iint\frac{f(x,y_1)-f(x,y_2)}{(1-{\rm i}\e y_1)(1-{\rm i}\e y_2)(1-{\rm i}y_2)}
\,dE_{A}(x)\,dE_{B_1}(y_1)\,dE_{B_2}(y_2)\\
=\lim_{\e\to0}\upint\!\!\!\iint\frac{f(x,y_1)}{(1-{\rm i}\e y_1)(1-{\rm i}\e y_2)(1-{\rm i}y_2)}
\,dE_{A}(x)\,dE_{B_1}(y_1)\,dE_{B_2}(y_2)\\
-\lim_{\e\to0}\upint\!\!\!\iint\frac{f(x,y_2)}{(1-{\rm i}\e y_1)(1-{\rm i}\e y_2)(1-{\rm i}y_2)}
\,dE_{A}(x)\,dE_{B_1}(y_1)\,dE_{B_2}(y_2)=X_\e-Y_\e,
\end{multline*}
where
$$
X_\e\df\lim_{\e\to0}\left(\iint\frac{f(x,y_1)}{1-{\rm i}\e y_1}
\,dE_{A}(x)\,dE_{B_1}(y_1)\right)(I-{\rm i}\e B_2)^{-1}(I-{\rm i}B_2)^{-1}
$$
and
$$
Y_\e\df\lim_{\e\to0}\left(\iint\frac{f(x,y_2)}{(1-{\rm i}\e y_2)(1-{\rm i}y_2)}
\,dE_{A}(x)(I-{\rm i}\e B_1)^{-1}\,dE_{B_2}(y_2)\right).
$$

%=\lim_{\e\to0}\left(\iint\frac{f(x,y_1)}{1-{\rm i}\e y_1}
%\,dE_{A}(x)\,dE_{B_1}(y_1)\right)(I-{\rm i}\e B_2)^{-1}(I-{\rm i}B_2)^{-1}\\
%-\lim_{\e\to0}\left(\iint\frac{f(x,y_2)}{(1-{\rm i}\e y_2)(1-{\rm i}y_2)}
%\,dE_{A}(x)(I-{\rm i}\e B_1)^{-1}\,dE_{B_2}(y_2)\right)\\

Clearly,
\begin{multline*}
\lim_{\e\to0} X_\e\\
=\lim_{\e\to0}\left(\iint\frac{f(x,y_1)}{1-{\rm i}y_1}
\,dE_{A}(x)\,dE_{B_1}(y_1)\right)(I-{\rm i}\e B_1)^{-1}(I-{\rm i}B_1)(I-{\rm i}B_2)^{-1}(I-{\rm i}\e B_2)^{-1}\\
=\left(\iint\frac{f(x,y_1)}{1-{\rm i}y_1}
\,dE_{A}(x)\,dE_{B_1}(y_1)\right)(I-{\rm i}B_1)(I-{\rm i}B_2)^{-1}\\
=f(A,B_1)(I-{\rm i}B_2)^{-1}
\end{multline*}
in the strong operator topology.
To compute the limit $\lim\limits_{\e\to0} Y_\e$, we observe that
\begin{align*}
\lim_{\e\to0}\iint\frac{f(x,y_2)}{1-{\rm i} y_2}
\,dE_{A}(x)(I-{\rm i}\e B_1)^{-1}\,dE_{B_2}(y_2)
=\iint\frac{f(x,y_2)}{1-{\rm i} y_2}
\,dE_{A}(x)\,dE_{B_2}(y_2)
\end{align*}
in the strong operator topology by Proposition 3.3 of \cite{BS3}. Hence,

\begin{align*}
\lim_{\e\to0} Y_\e=\lim_{\e\to0}\left(\iint\frac{f(x,y_2)}{1-{\rm i} y_2}
\,dE_{A}(x)(I-{\rm i}\e B_1)^{-1}\,dE_{B_2}(y_2)\right)(I-{\rm i}\e B_2)^{-1}\\
=\iint\frac{f(x,y_2)}{1-{\rm i} y_2}
\,dE_{A}(x)\,dE_{B_2}(y_2)=f(A,B_2)(1-{\rm i}B_2)^{-1}.
\end{align*}
in the strong operator topology. Thus, we have proved that

\begin{align*}
\lim_{\e\to0}\left(\upint\!\!\!\iint\frac{f(x,y_1)-f(x,y_2)}{y_1-y_2}
\,dE_{A}(x)\,dE_{B_1}(y_1)(B_1(\e)-B_2(\e))\,dE_{B_2}(y_2)\right)(I-{\rm i}B_2)^{-1}\\
=(f(A,B_1)-f(A,B_2))(I-{\rm i}B_2)^{-1}
\end{align*}
in the strong operator topology. $\bl$.

\begin{cor}
\label{cor510}
Let $f\in \E^\be_\s(\R^2)$, $1\le p\le2$. Suppose that $A$, $B_1$ and $B_2$ are self-adjoint operators such that 
$B_2-B_1\in\bS_p$. Then the following inequality holds:
$$
\|f(A,B_1)-f(A,B_2)\|_{\bS_p}\le\const\s\|f\|_{L^\be(\R^2)}\|B_1-B_2\|_{\bS_p}.
$$
\end{cor}

\Pf Applying \rf{eshchorazW} and \rf{otsenkaD2}, we obtain
\bey
\|f(A,B_1)-f(A,B_2)\|_{\bS_p}
\le\big\|\big(\dg^{[2]}f\big)(x,y_1,y_2)\big\|_{\B^\be\otimes_{\rm h}\B^\be\otimes^{\rm h}\B^\be}\|B_1-B_2\|_{\bS_p}\\
\le\const\s\|f\|_{L^\be(\R^2)}\|A_1-A_2\|_{\bS_p}.\quad\bl
\eey

\begin{cor} 
\label{511}
Let $f\in \Bs(\R^2)$, $1\le p\le2$. Suppose that $A$, $B_1$ and $B_2$ are self-adjoint operators such that 
$B_2-B_1\in\bS_p$. Then the following inequality holds:
$$
\|f(A,B_1)-f(A,B_2)\|_{\bS_p}\le\const\|f\|_{\Bs}\|B_1-B_2\|_{\bS_p}.
$$
\end{cor}

\Pf Let $f^{[0]}$ and $f_n$ be the functions defined by \rf{fnol'} and \rf{fn}.
It is easy to see that if $u\in D(B_1)=D(B_2)$, then
$$
f(A,B_1)u-f(A,B_2)u=f^{[0]}(A, B_1)u-f^{[0]}(A,B_2)u+\sum_{n=1}^\infty(f_n(A,B_1)u-f_n(A,B_2)u).
$$
It remains to observe that
$$
\|f^{[0]}(A,B_1)u-f^{[0]}(A,B_2)u\|\le\const\|f^{[0]}\|_{L^\infty}\|B_1-B_2\|\cdot\|u\|
$$
and
$$
\|f_n(A,B_1)u-f_n(A,B_2)u\|\le\const 2^n\|f_n\|_{L^\infty}\|B_1-B_2\|\cdot\|u\|
$$
by Corollary \ref{cor510}. $\bl$

\begin{thm}
\label{lipschitseva_otsenka_dlya_Besova}
Let $f\in\Bs(\R^2)$, $1\le p\le2$. Suppose that $A_1$, $A_2$, $B_1$ and $B_2$ are self-adjoint operators such that 
$A_2-A_1\in\bS_p$ and $B_2-B_1\in\bS_p$. Then the following inequality holds:
$$
\|f(A_1,B_1)-f(A_2,B_2)\|_{\bS_p}\le\const\|f\|_{\Bs}
\max\big\{\|A_1-A_2\|_{\bS_p},\|B_1-B_2\|_{\bS_p}\big\}.
$$
\end{thm}

\Pf We have
$$
\|f(A_1,B_1)-f(A_2,B_2)\|_{\bS_p}\le\|f(A_1,B_1)-f(A_2,B_1)\|_{\bS_p}+\|f(A_2,B_1)-f(A_2,B_2)\|_{\bS_p}.
$$
It remains to apply Corollary \ref{cor58}  and \ref{511}. $\bl$

\begin{thm} 
\label{ab12}
Let $f\in\Bs(\R^2)$. Suppose that $A_1$, $A_2$, $B_1$ and $B$ are self-adjoint operators such that 
$A_1-A_2\in\bS_2$ and $B_1-B_2\in\bS_2$. Then the following identity holds:
\begin{align}
\label{tryokhetazhka}
f(A_1,B_1)&-f(A_2,B_2)\nonumber\\[.2cm]
&=
\iint\!\!\upint\frac{f(x_1,y)-f(x_2,y)}{x_1-x_2}
\,dE_{A_1}(x_1)(A_1-A_2)\,dE_{A_2}(x_2)\,dE_{B_1}(y),\nonumber\\[.2cm]
&+\upint\!\!\!\iint\frac{f(x,y_1)-f(x,y_2)}{y_1-y_2}
\,dE_{A_2}(x)\,dE_{B_1}(y_1)(B_1-B_2)\,dE_{B_2}(y_2).
\end{align}
\end{thm}

\Pf Assume first that $f\in \E^\be_\s(\R^2)$. By Theorem \ref{teor1},
we have
\begin{align*}
f(A_1,B_1)&-f(A_2,B_1)\\[.2cm]
&=
\iint\!\!\upint\frac{f(x_1,y)-f(x_2,y)}{x_1-x_2}
\,dE_{A_1}(x_1)(A_1-A_2)\,dE_{A_2}(x_2)\,dE_{B_1}(y)
\end{align*}
while by Theorem \ref{teor2}, we have
\begin{align*}
f(A_2,B_1)&-f(A_2,B_2)\nonumber\\[.2cm]
&=
\upint\!\!\!\iint\frac{f(x,y_1)-f(x,y_2)}{y_1-y_2}
\,dE_{A_2}(x)\,dE_{B_1}(y_1)(B_1-B_2)\,dE_{B_2}(y_2).
\end{align*}
It remains to observe that
$$
f(A_1,B_1)-f(A_2,B_2)=(f(A_1,B_1)-f(A_2,B_1))+(f(A_2,B_1)-f(A_2,B_2)).
$$

Clearly, the set $\bigcup_{\s>0}\E^\be_\s(\R^2)$ is dense in $\Bs(\R^2)$.
It remains to observe that the set of functions $f$, for which \rf{tryokhetazhka} holds is closed in $\Bs(\R^2)$.
This follows from Theorem \ref{lipschitseva_otsenka_dlya_Besova} and estimates \rf{otsenkaD3} and  \rf{otsenkaD4}. $\bl$

In the same way we can obtain the following result:

\begin{thm} 
\label{ba21}
Let $f\in\Bs(\R^2)$. Suppose that $A_1$, $A_2$, $B_1$ and $B$ are self-adjoint operators such that 
$A_1-A_2\in\bS_2$ and $B_1-B_2\in\bS_2$. Then the following identity holds:
\begin{align*}
f(A_1,B_1)&-f(A_2,B_2)\nonumber\\[.2cm]
&=
\iint\!\!\upint\frac{f(x_1,y)-f(x_2,y)}{x_1-x_2}
\,dE_{A_1}(x_1)(A_1-A_2)\,dE_{A_2}(x_2)\,dE_{B_2}(y),\nonumber\\[.2cm]
&+\upint\!\!\!\iint\frac{f(x,y_1)-f(x,y_2)}{y_1-y_2}
\,dE_{A_1}(x)\,dE_{B_1}(y_1)(B_1-B_2)\,dE_{B_2}(y_2).
\end{align*}
\end{thm}

\

\
 
 \begin{footnotesize}
 
\noindent
\begin{tabular}{p{7cm}p{15cm}}
A.B. Aleksandrov & V.V. Peller \\
St.Petersburg Department & St.Petersburg State University\\
Steklov Institute of Mathematics  & Universitetskaya nab., 7/9\\
Fontanka 27, 191023 St.Petersburg & 199034 St.Petersburg \\
Russia&Russia\\

\\

 St.Petersburg State University&Department of Mathematics\\
 Universitetskaya nab., 7/9&Michigan State University\\
199034 St.Petersburg&East Lansing, Michigan 48824\\
Russia&USA\\
email: alex@pdmi.ras.ru&\\
&St.Petersburg Depatment\\
&Steklov Institute of Mathematics\\
&Fontanka 27, 191023 St.Petersburg\\
&Russia\\
&and\\
&Peoples' Friendship University\\
& of Russia (RUDN University)\\
&6 Miklukho-Maklaya St., Moscow,\\
& 117198, Russian Federation\\
& email: peller@math.msu.edu
\end{tabular}

\end{footnotesize}

\end{document}